\documentclass[12pt]{amsart}
\usepackage{amscd,amssymb,verbatim}
\usepackage{amssymb,amsmath,amsthm,epsfig}

\def\cnum#1{\bigcirc\kern -8pt#1}
\newcommand{\n}{\noindent}
\newcommand{\vp}{\varepsilon}
\newcommand{\ovl}{\overline}

\theoremstyle{plain}
\newtheorem{thm}{Theorem}[section]
\newtheorem{pro}[thm]{Proposition}
\newtheorem{lem}[thm]{Lemma}

\theoremstyle{definition}
\newtheorem*{defn}{Definition}

\theoremstyle{remark}
\newtheorem*{rem}{Remark}

\newcommand{\e}{\varepsilon}

\newcommand{\supp}{{\rm supp}\,}

\newcommand{\aco}{{\rm aco}\,}
\newcommand{\N}{\mathbb{N}}

\newcommand{\R}{\mathbb{R}}

\begin{document}

\title{Strictly singular, non-compact operators exist on the   space of
        Gowers and Maurey}

\author{G.  Androulakis \and Th. Schlumprecht}

\
\thanks{Research supported by NSF}


\date{1/22/2001}

\maketitle 
 \markboth{G.  Androulakis \and Th. Schlumprecht}{Strictly singular, non-compact operators on Gowers' and Maurey's   space}

\noindent
{\bf Abstract}
We construct a strictly singular non-compact operator on 
 Gowers' and  Maurey's space $GM$. 

\section{Introduction} \label{sec1}

\indent

In 1993 W.T.~Gowers and B.~Maurey \cite{GM}  solved
 the famous ``unconditional basic sequence problem'' by constructing the first 
known example of a
space that does not contain any unconditional basic sequence. In
the present paper this space is denoted by $GM$. Furthermore it was shown in 
\cite{GM}
 that the space $GM$ is {\em hereditarily
indecomposable} (HI), i.e.\ no infinite dimensional subspace of $GM$ can be
decomposed into a direct sum of two further infinite dimensional closed 
subspaces.
 As shown in \cite{GM} every bounded operator on a complex HI space can be 
written as
 a sum of a  multiple of the identity  and a {\em strictly compact operator}. 
Recall that an operator
$T$ is strictly singular if no restriction of $T$ to an infinite
dimensional subspace is an isomorphism.
 Actually Lemma 22 of \cite{GM} implies immediately that the real version of 
$GM$ has
 also the property that every operator on a subspace of it is a strictly singular perturbation
of  a multiple of the identity.
 In 
\cite{GM}  it is asked whether or not every
operator on $GM$ can be written as a compact perturbation of a multiple of the
identity. If the answer to this question were positive then by \cite{AS} the
space $GM$ would be the first known example of an infinite dimensional Banach
space such that every operator on it has a non-trivial invariant subspace
i.e.\ $GM$ would be a positive solution to the invariant subspace problem. 

In 1999 W.T.~Gowers showed that strictly singular non-compact operators can be
defined on certain subspaces of $GM$, \cite{G}. Our main
result is

\begin{thm}\label{thm1.1}
There exists a strictly singular non-compact operator $T\colon \ GM\to GM$.
\end{thm}

Since strictly singular non compact operators on an infinite dimensional banach space cannot be written as
compact perturbations of a multiple of the identity, we give a
negative  answer to  the question  of W.T.~Gowers
and B.~Maurey. Also, we show that the space of operators on $GM$ 
contains a subspace isometric to $\ell_\infty$, the
Banach space of all bounded sequences of scalars. Thus the set of operators
on $GM$ is non-separable in contrast to the separability of the  set
of compact perturbations of a multiple of the identity on a space with a basis. 

Concerning the invariant subspace problem,  note that if every strictly singular 
operator
on $GM$ has some compact power, then by \cite{L}, $GM$ would still be a positive
solution to the invariant subspace problem. However, we conjecture that there
exists a strictly singular operator  on $GM$ none of whose powers is compact.

The construction of the space $GM$ is based on the space $S$, the first
known example of an arbitrarily distortable space, constructed by the second
named author in \cite{S1}.

The space $S$ was used by E.~Odell and the second  author in order to
show  that the separable Hilbert space  is arbitrarily
distortable \cite{OS}. The second named author proved in \cite{S2} that the
space $S$ is complementably minimal.
This means that for every infinite dimensional subspace of $S$ there exists a
further subspace which is isomorphic to $S$ and complemented in $S$. Note
that if a Banach space $X$ with the Aproximation Property (AP) is
complementably  minimal then every non-trivial
closed two-sided operator ideal $\mathcal{I}$ of $X$ satisfies $\mathcal{K}
\subseteq \mathcal{I} \subseteq \mathcal{S}$ where $\mathcal{K}$ is the ideal
of compact operators and $\mathcal{S}$ is the ideal of strictly singular
operators on $X$. Indeed, since $\mathcal{I}$ is closed and
non-trivial, and $X$ has the AP we
have $\mathcal{K}\subseteq\mathcal{I}$. Assume that there exists $T\in
\mathcal{I}\backslash \mathcal{S}$. Since $T$ is not strictly singular, there
exists an infinite dimensional subspace $Y$ of $X$ such that $T$ restricted on
$Y$ is an isomorphism. Since $X$ is complementably minimal there exists a
subspace $Z$ of $T(Y)$ such that $Z$ is isomorphic to $X$ and complemented in
$X$. Let $j\colon \ X\to Z$ be an onto isomorphism and let $P\colon \ X\to Z$
be an onto projection. Finally note that $ATB$ is the identity  on $X$
where $A  = j^{-1}P$ and $B = (T|_Z)^{-1}j$ are operators on $X$. Therefore
$\text{id } \in \mathcal{I}$
and $\mathcal{I}$ consists of all operators on $X$. Hence, for a complementably minimal space $X$
having (AP),
the lack of strictly singular non-compact operators on $X$ is equivalent to
$X$ being {\em simple}. Recall that a Banach space $X$ is simple if the only
two-sided closed operator ideal on $X$ is the ideal of compact operators. The
only known simple spaces are the spaces $\ell_p$ $(1\le p<\infty)$ and $c_0$
\cite{H} (see also \cite{FGM} and \cite{P} page 82). W.B.~Johnson
asked  whether or not $S$ is simple. We prove

\begin{thm}\label{thm1.2}
There exists a strictly singular non-compact operator on $S$.
\end{thm}

 Thus $S$ is not
a simple space. Also Theorem~\ref{thm1.1} implies that $GM$ is not a simple
space. We do not know how many  two-sided  closed operator ideals
exist on $S$ and on $GM$.

It will follow from our work that formaly the same operator $T$ can be 
considered
in both Theorems~\ref{thm1.1} and \ref{thm1.2} (either as an operator on $GM$
or as an operator on $S$). The operator $T$ has the form $T = \Sigma
x^*_i\otimes e_i$ where $(x^*_i)$ is a seminormalized block sequence  in $S^*$ 
as well as in $GM^*$ and $(e_i)$ is
the unit  vector basis of $S$ or $GM$ respectively. 

We assume all our Banach spaces being defined over $\R$ noting 
that the result can be easily transfered to the complex case.

We now recall the definition of $S$ and will introduce first some basic 
notation.
 $c_{00}$ is the vectorspace of sequences in $\R$ for which only finitely
 many coordinates are not zero. 
  For $x\in c_{00}$ the {\em support} of $x$ is the set $\{i\in\N: x_i\not=0\}$ 
which
  we denote by supp$(x)$.
 $(e_i)$ is the usual unit basis in
 $c_{00}$, i.e. $e_i(j)=\delta_{ij}$, for $i,j\in\N$.
 If $x=\sum x_ie_i\in c_{00}$ and $E\subset \N$ we write
  $E(x)=\sum_{i\in E} e_i$. For two finite sets
$E,F\subset \N$
  we write $E<F$ if $\max E<\min F$ ($\max \emptyset=0$), and for $x,y\in c_{00}$ we write $x<y$ if
  $\text{supp}(x)< \text{supp}(y)$.

Let $\|\cdot\|_{\ell_p}$ denote the usual norm on $\ell_p$ if $1\le p\le 
\infty$.
 Let $f$ denote the function $f(x) = \log_2(x+1)$.  $S$ is the completion
 of $(c_{00},\|\cdot\|_S)$ and  $\|\cdot\|_S$ is the
unique norm on  $c_{00}$ which satisfies the implicit formula:
$$\|x\|_S = \|x\|_{\ell_\infty} \vee \sup_{\stackrel{\scriptstyle 2\le n\in
\N}{\scriptstyle E_1 < E_2 <\cdots< E_n}} \frac1{f(n)} \sum^n_{j=1}
\|E_jx\|_S$$
 The unique existence of such a norm is easy to show and it is also
 not hard to prove that $(e_i)$ is a 1-unconditional
  and 1-subsymmetric basis of $S$ (see \cite{S1}).

 For $\ell\in\N$, $\ell\ge 2$ we define the equivalent norm $\|\cdot\|_\ell$ on $S$
by 
$$\|x\|_\ell = \sup_{E_1 <\cdots< E_\ell} \frac1{f(\ell)} 
\sum^\ell_{j=1}\|E_jX\|_S.$$
The construction of the space $GM$ will be recalled in Section~3.

\section{Existence of strictly singular,  non-compact operators on
$\pmb{S}$}\label{sec2}

The main goal of this section is to prove Theorem \ref{thm1.2}. We
start by stating a sufficient condition for an operator $T\colon \ S\to S$, of
the form
$$T = \sum_i x^*_i\otimes e_i,$$
to be strictly singular but not compact. 

\begin{pro}\label{pro2.2}
Assume that $(x^*_i)^\infty_{i=1}$ is a seminormalized block sequence in $S^*$
and that there is an increasing sequence $(C(\ell))_{\ell\in\N}\subset \R^+$, with $C(\ell) 
\nearrow \infty$, if $\ell\nearrow \infty$, for
which the following condition holds.
\begin{align} \label{eq2.1}
& \mbox{If $(z_i)^\infty_{i=1}$ is a block sequence in $S$, so that for
each $i\in \N$, $x^*_i(z_i)=1$ and $x^*_{i-1}<z_i$}\\
& <x^*_{i+1} \mbox{ (take $x^*_0 = 0$),
 then for any $2\le \ell \in \N$ and 
$(\lambda_i)^\infty_{i=1}\in c_{00}$ we have that } \nonumber \\
& \hskip 1.8in \left\|\sum^\infty_{i=1} \lambda_i e_i\right\|_\ell \le 
\frac1{C(\ell)}
\left\|\sum^\infty_{i=1} \lambda_iz_i\right\|_S. \nonumber
\end{align}
Then the operator $T = \sum^\infty_{i=1} x^*_i \otimes e_i$, with
$T(x) = \sum^\infty_{i=1} x^*_i(x) \cdot e_i$, for $x\in S$, is
bounded, strictly singular, but not compact. 
\end{pro}

\begin{proof}
In order to see that $T\colon \ S\to S$ is bounded let $x\in c_{00}$.
 We can write $x$ as $x = \sum^\infty_{i=1}
\lambda_iz_i$, where $(z_i)$ is a block sequence in $S$ so that $x^*_i(z_i)=1$
and $x^*_{i-1}<z_i<x^*_{i+1}$. Thus $Tx = \Sigma \lambda_i
e_i$. If $\|Tx\|_S = \|Tx\|_{\ell_\infty}$, then
$$
\|Tx\|_{\ell_\infty} = \max_{i\in \N} |\lambda_i|
= \max_{i\in \N} \left|x^*_i\left(\sum \lambda_iz_j\right)\right|
\le \sup_{i\in \N} \|x^*_i\|_{S^*} \|x\|_S.
$$
If $\|Tx\|_S = \|Tx\|_\ell$ for some $\ell\in \N$, $\ell\ge 2$, then
it follows from our assumption (\ref{eq2.1}) that
$$\|Tx\|_\ell = \left\|\sum^\infty_{i=1} \lambda_i e_i\right\|_\ell \le
\frac1{C(\ell)} \left\|\sum \lambda_iz_i\right\|_S \le \frac1{C(2)} \|x\|_S.$$
In order to show that $T$ is strictly singular we consider an arbitrary
infinite dimensional subspace $X$ of $S$, and let $\ell_0\in \N$. $X$
contains an element $x$ for which there exists an $\ell\ge \ell_0$ so that
$\|x\|_S = \|x\|_\ell$. Indeed, as it was shown in \cite{S1}, we find for any
$N\in \N$ and $\vp>0$ a normalized block $(y_i)^N_{i=1}$ in $X$ which
is $(1+\vp)$-equivalent to the $\ell^N_1$-unit vector basis, in particular it
follows that
$$\frac1{1+\vp} \le \left\|\frac1N \sum^N_{i=1} y_i\right\|_S.$$
On the other hand it was shown in \cite{S1} that given $\ell_0\in \N$
we can choose $N$ big enough so that for any $m\in \N$, $m\le \ell_0$, it 
follows that
$\|\frac1N \sum^N_{i=1} y_i\|_m <
(1+\delta(\vp))/f(m)$, with $\delta(\vp) \downarrow 0$, for
$\vp\downarrow 0$. Since clearly 
$\|\frac1N \sum^N_{i=1} y_i\|_{\ell_\infty} \le \frac1N$ there  exists an 
$\ell\ge \ell_0$ so that
$\|\frac1N \sum^N_{i=1} y_i\|_\ell = \|\frac1N  \sum^N_{i=1} y_i\|_S$. 
Now assume that $X$ is an infinite
dimensional subspace of $S$ on which $T$ acts as an isomorphism. For any
$\ell_0$ we can choose a $y\in T(X)$, $y\in T(x)$, so that $\|y\| =
\|y\|_\ell$, with $\ell\ge\ell_0$. As before we write $x$ as $x =
\sum^\infty_{i=1} \lambda_iz_i$ where $(z_i)$ is a block sequence with
$x^*_i(z_i)=1$ and $x^*_{i-1}<z_i<x^*_{i+1}$ for $i\in \N$. 
Then
$\|Tx\|_S = \|Tx\|_\ell = \left\|\sum^\infty_{i=1} \lambda_i e_i
\right\|_\ell \le$ $\frac1{C(\ell)} \left\|\sum^\infty_{i=1} \lambda_i
z_i\right\|_S \le \frac1{C(\ell_0)} \|x\|_S$. 
Since $\ell_0\in \N$ was arbitrary and since $C(\ell) \nearrow \infty$
for $\ell\nearrow \infty$ we arrive to a contradiction to the assumption that
$T|_X$ in an isomorphism. 

Finally, $T$ is not compact since $(x^*_i)$ is a
seminormalized block sequence.
\end{proof}

\begin{rem}
Using the same proof, Proposition~\ref{pro2.2} can be generalized as follows: 
Assume $(x^*_i)^\infty_{i=1} \subset S^*$ and $(y_i)^\infty_{i=1}
\subset S$ are seminormalized block sequences and $(C(\ell))\subset\R_+$ 
so that $C(\ell)\nearrow \infty$ if $\ell\nearrow \infty$ and
the following conditions holds:\newline
If $(z_i)^\infty_{i=1}$ is a block sequence in $S$, so that for each $i\in
\N$, $x^*_i(z_i) = 1$ and 
$x^*_{i-1}<z_i<x^*_{i+1}$,
then for any $2\le \ell \in \N$ and $(\lambda_i)^\infty_{i=1}\in
c_{00}$ we have that
$$\left\|\sum \lambda_iy_i\right\|_\ell \le \frac1{C(\ell)} \left\|\sum
\lambda_iz_i\right\|_S   + \max|\lambda_i|.$$
Then the operator $T = \Sigma x^*_i\otimes y_i$ is bounded, strictly singular
but not compact.
\end{rem}

\begin{rem}  The rest of this section is devoted to the construction of 
 a semi normalized block $(x^*_n)$ in $S^*$  and a sequence $C(\ell)_{\ell\in\N}$ satisfying the condition stated in
  Proposition~\ref{pro2.2}. To do that we will have to overcome several 
technical 
   difficulties. In a first reading one can proceed directly  to Section 3.
There we will
   show that by spreading out the coordinates  of the sequence $(x^*_i)$
   and  under some mild growth condition on $C(\ell)$ (which will be verified) the operator
   $\Sigma x^*_i\otimes y_i$ is well defined on $GM$, strictly singular and 
non-compact.
  \end{rem}  
In order to define the sequence $(x^*_m)\subset S^*$ 
 we will need the following notion of
{\em finitely branching trees.} Consider the set $\cup^\infty_{k=0}
\N^k$, the set of all finite sequences with values in $\N$,
which will be partially ordered by extensions, i.e.\ if $\mu = 
(\mu_1,\ldots,\mu_m)$ 
and $\nu = (\nu_1,\ldots, \nu_n)$ are in $\cup^\infty_{k=0}
\N^k$ we write $\mu \prec\nu$ if $m\le n$ and $\mu_1=\nu_1$,
$\mu_2=\nu_2,\ldots, \mu_m = \nu_m$. For $\mu = (\mu_1,\ldots, \mu_m)\in
\cup^\infty_{k=0} \N^k$, we call $m$ the {\em length of\/}
$\mu$ and denote it by $|\mu|$.

A non-empty subset $\mathcal{T}$ of $\cup^\infty_{k=0} \N^k$
is called a {\em finitely branching tree}, or simply a tree (since we will not
deal with other kinds of trees) if the following two conditions hold. 
\begin{equation}
 \text{If\ } \mu = (\mu_1,\ldots, \mu_m)\in \mathcal{T}\text{\ and\ } 0 \le k 
\le m,
 \text{\ then\ also\ }(\mu_1,\ldots,\mu_k) \text{\ lies in\ }
 {\mathcal T}.
\end{equation}
Since $\mathcal{T}$ is non-empty, this implies in particular that $\emptyset
\in \mathcal{T}$.
\begin{align}
& \mbox{If $\mu\in \mathcal{T}$, then either $\mu$
is maximal in $\mathcal{T}$, i.e.\ no extension of $\mu$ lies in  $\mathcal{T}$, 
or there is a }\\ 
& k_\mu=k_\mu^{\mathcal T}\in \N \text{ so that } 
 \{(\mu,i)\colon \ i\in \{1,2,\ldots, k_\mu\}\} = 
\{\nu\in \mathcal{T}\colon \ \mu\prec \nu, |\nu| = |\mu|+1\}. \nonumber
\end{align}

We call the elements of trees {\em nodes}. If $\mathcal{T}$ is a tree and $\mu
= (\mu_1,\ldots, \mu_{|\mu|})\in \mathcal{T}\backslash \{\emptyset\}$ we call
$(\mu_1,\ldots, \mu_{|\mu|-1})$ the {\em immediate predecessor\/} of $\mu$. If
$\mu\in \mathcal{T}$ is not maximal in $\mathcal{T}$ we call the nodes
$(\mu,i)$, $i\le k_\mu^{\mathcal T}$, {\em immediate successors\/} of $\mu$ in
$\mathcal{T}$. 
Since we will always deal with only one tree $\mathcal T$ at a time we denote the 
numbers
 of successors of an non-maximal $\mu\in\mathcal T$ simply by $k_\mu$.
We say that a tree $\mathcal{T}$ is of length $\ell$, $\ell\in
\N_0$, if $\mathcal{T}\subset \cup^\ell_{k=0}
\N^k$ and all maximal nodes in $\mathcal{T}$ have length $\ell$. A tree
is said to have infinite length if it does not contain maximal nodes. On a
tree $\mathcal{T}$ having finite or infinite length we also introduce the {\em
lexicographic order\/} which we denote by $\prec_{\text{lex}}$, i.e.\ we
well-order $\mathcal{T}$ into: \ $\emptyset, (1), (2),\ldots, (k_\emptyset)$,
$(1,1),\ldots, (1,k_{(1)})$, $(2,1),\ldots, (2,k_{(2)})$, $(3,1),\ldots$~.
To a given tree $\mathcal{T}$ having finite length we want to define {\em
associated vectors\/} in $S$ and $S^*$. These vectors are defined up to {\em
equality in distribution\/}. For $x,y\in c_{00}$ we say that $x$ and $y$ have
the same distribution and write $x =_{\text{dist}} y$ if for some $k\in
\N$, $(\alpha_i)^k_{i=1} \subset \R$, and 
$m_1 < m_2<\ldots<m_k$ and $n_1 < n_2 <\cdots < n_k$ in $\N$ it follows that
 $x = \sum^k_{i=1}\alpha_i e_{m_i}\quad \text{and}\quad y = \sum^k_{i=1}
\alpha_i e_{n_i}$.

\begin{defn}
Let $\mathcal{T}$ be a tree of length $\ell\in \N_0$. For
$k\le\ell$ and each $\mu = (\mu_1,\ldots, \mu_k)\in \mathcal{T}$ we define the
numbers
\begin{equation}
\label{eq2.4}
\alpha(\mu) = \prod^{k-1}_{i=0} \frac{f(k_{(\mu_1,\ldots, \mu_i)})}{k_{(\mu_1,
\ldots, \mu_i)}}
\mbox{\  and\ \  }
\beta(\mu) = \prod^{k-1}_{i=0} \frac1{f(k_{(\mu_1,\ldots,\mu_i)})}
\end{equation}
(if $k\!=\!0$, then $(\mu_1,\ldots, \mu_0)\! =\! \emptyset$ and 
$\alpha(\emptyset)\! =\!
\beta(\emptyset)\! = \!1$). We say that $x\in S$ is {\em associated to\/}
$\mathcal{T}$ if
\begin{equation}\label{eq2.6}
x =\sum_{\mu\in \mathcal{T}, |\mu|=\ell} \alpha(\mu)\cdot e_{n(\mu)}
\end{equation}
and we say that $x^*\in S^*$ is {\em associated to\/} $\mathcal{T}$ if
\begin{equation}\label{eq2.7}
x^* = \sum_{\mu\in \mathcal{T}, |\mu| =\ell} \beta(\mu) e^*_{n(\mu)}
\end{equation}
where $(n(\mu))_{\mu\in\mathcal T}\in \N$ has the property that $n(\mu) < n(\nu)$ if
$\mu \prec_{\text{lex}} \nu$.
\end{defn}

\begin{rem}
In the above notation it is easy to see that $\|x\|_S \ge 1$, $\|x^*\|_{S^*} \le 
1$
and that $x^*(x) = 1$ (if the $n(\mu)$'s coincide).
\end{rem}

\begin{rem}
There is also a recursive way to introduce the vectors $x$ and $x^*$ which are
associated to a tree of length $\ell$. If $\ell=0$, i.e.\ $\mathcal{T} =
\{\emptyset\}$, the associated vectors in $S$ and $S^*$ are simply the 
elements of the unit vector basis. 
Assume that for an $\ell\ge 0$, we have defined the vectors associated to
trees of length $\ell$ in $S$ and $S^*$ respectively. Let $\mathcal{T}$ be a
tree of length $\ell+1$.  For $i\le k_\emptyset$, the number of immediate
successors of $\emptyset$ in $\mathcal{T}$, we let
$\mathcal{T}_i = \left\{\mu \in \cup^\ell_{j=0} \N^j\mid (i,\mu)
\in \mathcal{T}\right\}$.
Then $\mathcal{T}_i$ is a tree of length $\ell$ and we choose vectors $x_i\in
S$ and $x^*_i\in S$ associated to $\mathcal{T}_i$. Furthermore we require
$(x_i)^{k_\emptyset}_{i=1}$ and $(x^*_i)^{k_\emptyset}_{i=1}$ to be blocks. 
Then the
vectors associated to $\mathcal{T}$ are
\begin{equation}
\label{eq2.8}
x = \frac{f(k_\emptyset)}{k_\emptyset} \sum^{k_\emptyset}_{i=1}
x_i\quad \text{and} \quad
x^* = \frac1{f(k_\emptyset)} \sum^{k_\emptyset}_{i=1} x^*_i.
\end{equation}
More generally if $0 \le k\le \ell$, and if $\mathcal{T}$ is a tree of length
$\ell$, then every $x\in S$, respectively $x^*\in S$, associated to
$\mathcal{T}$ can be written as
\begin{equation}
\label{eq2.10}
x = \sum_{\mu\in\mathcal{T}, |\mu|=k} \alpha(\mu) x_\mu,\quad
\text{and} \quad
x^* = \sum_{\mu\in \mathcal{T}, |\mu|=k} \beta(\mu) x^*_\mu,
\end{equation}
where $x_\mu$ and $x^*_\mu$ are associated to $\mathcal{T}_\mu := \left\{\nu
\in \cup^{\ell-k}_{j=0} \N^j \mid (\mu,\nu)\in
\mathcal{T}\right\}$ and $(x_\mu)$ and $(x^*_\mu)$ are blocks with respect to
the lexicographic order.
\end{rem}

Every tree $\mathcal{T}$ is determined, by the numbers $k_\mu$, for $\mu\in
\mathcal{T}$ being non-maximal. In order to choose the wanted block sequence
$(x^*_n) \subset S^*$ which satisfies the requirements of
Proposition~\ref{pro2.2} we will first choose a ``lacunary enough'' subset
$K\subset \N$. Secondly, we will choose an infinite tree $\mathcal{T}$,
with $(k_\mu)_{\mu\in\mathcal{T}} \subset K$ and such that $(k_\mu)$ increases
with respect to the lexicographic order. Then we will choose a sequence
$(L_n)^\infty_{n=1} \subset \N$, which increases ``fast enough'' to
$\infty$. Finally we let $(x^*_n)\subset S^*$ be a block sequence for which
$x^*_n$, $n\in \N$, is associated to $\mathcal{T}_{L_n} = \{\mu \in
\mathcal{T}\mid |\mu|\le L_n\}$. Let us formulate these conditions precisely.

Let $(\vp_i)^\infty_{i=1}$ be a decreasing sequence of positive numbers so
that $\sum^\infty_{i=1} \vp_i<\infty$. Assume that $(k_i)^\infty_{i=1}
\subset \N$ satisfies the following growth conditions
\begin{align}
\label{eq2.12}
&\frac23 f\left(\frac34 f(k_1)\right) \ge 1,\text{\ and\ } \ln k_1\ge 3\\
\label{eq2.13}
&\text{for}\quad r>k_1\quad \text{and}\quad a>1, \text{ then } f(r^a) \ge
\frac34 af(r).
\end{align}
Note that always $f(r^a) \le af(r)$, and that we can assume (\ref{eq2.13})
for a  large enough $k_1$ since $f(\cdot)$ is asymptotically logarithmic. For
$j\in \N$ we require the following inequalities
\begin{equation}
\label{eq2.14}
\left(\frac6{3f(k_j)} \prod^{j-1}_{s=1} \frac{k_s}{f(k_s)}\right)^{1/2}
\sum^\infty_{\ell=0} \left(\frac83 k_1\right)^{-\ell/2} 
 +
2 \frac{\log_2 f(k_j)}{f(k_j)^{1/2}} \left(\prod^{j-1}_{s=1}
\frac{k_s}{f(k_s)}\right)^{1/2} < \vp_j.
\end{equation}
For $r>1$ we define
\begin{equation}\label{eq2.16}
G(r) := \max_{i\in\N}
\frac{f(r) f(k_i)}{f(r k_i)} \prod^{i-1}_{j=1} \frac{f(k_j)}{k_j}
\end{equation}
and establish several properties. First note that above maximum is welldefined 
since 
 for fixed $r>0$ 
 $\frac{f(r)f(k_i)}{f(r k_i)} \prod^{i-1}_{j=1} \frac{f(k_j)}{k_j}$
  converges to zero if $i\nearrow\infty$.

\begin{pro}\label{pro2.3} 
\phantom{aaaaaaaaaaaa} \hfil\break
\vspace{-15pt}
\begin{itemize}
\item[a)] $G(r)$ is increasing in $r>1$.
\item[b)] If $(k_{i_j})^\infty_{j=1}$ is a subsequence of $(k_i)^\infty_{i=1}$
and we define
$$\ovl G(r) = \max_{j\in \N} \frac{f(r) f(k_{i_j})}{f(rk_{i_j})}
     \prod^{j-1}_{s=1} \frac{f(k_{i_s})}{k_{i_s}}$$
(i.e.\ we replace $(k_i)$ in the definition of $G$ by $(k_{i_j})$), then $\ovl
G(r) \ge G(r)$.
\end{itemize}
\end{pro}

\begin{proof}
For (a) we have to show that the function 
$[1,\infty)\ni x\mapsto\frac{f(x)}{f(ax)}$,
 $a\ge k_1$ is increasing. By
taking derivatives we need to show that
$$\frac{\ln(1+ax)}{1+x}-\frac{a\ln(1+x)}{1+ax}\ge 0$$.
 
An easy  computation shows that:
\begin{equation*}
 \frac{\ln(1+ax)}{1+x}-\frac{a\ln(1+x)}{1+ax}\ge \frac{\ln (ax)}{1+x}-\frac{a\ln(x)+ a}{1+ax}
  \ge\frac{ax(\ln(a)-3)}{(1+x)(1+ax)}\ge 0\text{ (by (\ref{eq2.12})).}
 \end{equation*}

In order to show (b) let $r>1$ and choose $i\in \N$ such that
$$G(r) =  \frac{f(r) f(k_i)}{f(rk_i)}\prod^{i-1}_{s=1} \frac{f(k_s)}{k_s}.$$
Secondly choose an $s\ge 1$ so that
 $k_{i_{s-1}} \le k_{i-1} < k_i \le k_{i_s}$
$(k_0 = k _{i_0} = 1)$.
Then
$$
\ovl G(r) \ge  \frac{f(r)f(k_{i_s})}{f(rk_{i_s})}\prod^{s-1}_{t=1} 
\frac{f(k_{i_t})}{k_{i_t}}
\ge  \frac{f(r) f(k_{i_s})}{f(r k_{i_s})}\prod^{i-1}_{t=1} \frac{f(k_t)}{k_t}
\ge \frac{f(r)f(k_i)}{f(rk_i)}\prod^{i-1}_{t=1} \frac{f(k_t)}{k_t}  = G(r).
$$
\end{proof}

\begin{lem}\label{lem2.4}
For $i\in\N$ choose $m_i\in \R^+$ such that $f(m_ik_i) = k_i$,
i.e.\ $m_i = \frac{2^{k_i}-1}{k_i}$.
Assume that $r\in [m_{i-1}, m_i)$, for some $i\ge 2$ and define the sequence
$(r_\ell)^\infty_{\ell=0}$ inductively by $r_0=r$ and $r_{\ell+1} =
r^{f(r_\ell)}_\ell$. Then
\begin{equation}\label{eq2.17}
\sum^\infty_{\ell=0} \frac1{\sqrt{G(r_\ell)}} \le \sum^\infty_{j=i} (\vp_j +
\vp_{j-1}).
\end{equation}
\end{lem}

The reader interested in the technical details of the proof is referred to
 the end of this section.

In addition to the growth conditions on $(k_i)^\infty_{i=1}$ stated in
(\ref{eq2.12}), (\ref{eq2.13}), and (\ref{eq2.14}) we will need
a further property. Recall from \cite{S1} that for given $\ell\in \N$
and $\vp>0$ one can choose $n_1 < n_2<\ldots < n_\ell$ fast enough increasing
so that a block $(y_i)^\ell_{i=1}$, with $y_i =_{\text{dist}}
\frac{f(n_i)}{n_i} \sum^{n_i}_{j=1} e_j$, for $i=1,\ldots, \ell$, is
$(1+\vp)$-equivalent to $(e_i)^\ell_{i=1}$ in $S$. In particular this
implies that $\left\|\sum^\ell_{i=1} y_i\right\|_S \le
\frac{(1+\vp)\ell}{f(\ell)}$. Therefore we can and will require the following
additional property on the sequence $(k_i)^\infty_{i=1}$.
\begin{align} \label{eq2.18}
&\mbox{If $\mathcal{T}$ is a tree of infinite length with $\{k_\mu\colon
 \mu\in\mathcal{T}\} \subset (k_i)^\infty_{i=1}$ and $k_\mu <\! k_\nu$, if 
$\mu\!\prec_{\text{lex}}\! \nu$, } \\ 
&\mbox{then for  $\ell\in \N$
and   $x\! \in\!  S$ associated to $\{\mu\! \in\! \mathcal{T}\colon \
|\mu|\le\ell\}$  it follows that $\|x\|\le 2$.} \nonumber
\end{align}
 Note that (\ref{eq2.18}) implies that for a tree $\mathcal{T}$
as in (\ref{eq2.18}), and $\ell\in \N$, any $x^*\in S^*$ associated to
$\{\mu\in \mathcal{T}, |\mu|\le\ell\}$ is at least of norm $\frac12$.
Assume now that $(k_i)^\infty_{i=1}$ satisfies (\ref{eq2.12}) through
(\ref{eq2.14}) and (\ref{eq2.18}), and let $G(r)$, $r>1$, be defined as in
(\ref{eq2.16}). We choose a sequence $(L_n)^\infty_{n=1}\subset \N$
which has the following property (\ref{eq2.19}).
\begin{align}\label{eq2.19}
& \mbox{ For $r>1$, the number } L_{\lfloor f(f(r)) \rfloor } \mbox
{ is big  enough to
ensure  that } \qquad \qquad \qquad \qquad \qquad \\
& \hskip 1.45in G(r) = \max_{i \le L_{\lfloor f(f(r)) \rfloor }}
\frac{f(k_i)f(r)}{f(k_ir)}\prod^{i-1}_{j=1} \frac{f(k_j)}{k_j}. \nonumber
\end{align}
Note that (\ref{eq2.19}) passes through to any subsequence
$(k_{i_s})^\infty_{s=1}$ of $(k_i)$, i.e.\
\begin{equation*}
\ovl G(r) := \max_{s\in\N} 
\frac{f(k_{i_s})f(r)}{f(k_{i_s}r)}\prod^{s-1}_{t=1}\frac{f(k_{i_t})}{k_{i_t}}
 = \max_{s\le L_{\lfloor f(f(r)) \rfloor }}
\frac{f(k_{i_s})f(r)}{f(k_{i_s}r)} \prod^{s-1}_{t=1} \frac{f(k_{i_t})}{k_{i_t}}.
\end{equation*}

Finally we choose an infinite tree $\mathcal{T}$ so that $\{k_\mu\colon \
\mu\in \mathcal{T}\}\subset (k_i)^\infty_{i=1}$ and so that $k_\mu < k_\nu$ if
$\mu \prec_{\text{lex}} \nu$.
 We let $(x^*_n)^\infty_{n=1}$ be a block
sequence in $S^*$, with $x^*_n$ being associated to the tree $\{\mu\in
\mathcal{T}\colon \ |\mu|\le L_n+1\}$.

\begin{lem}\label{lem2.5}
Assume that $(z_n)^\infty_{n=1}\subset S$ has the property that $x^*_n(z_n) =
1$, and $x^*_{n-1}<z_i<x^*_{n+1}$.
If $m\in\N$ and $I\subset \N$ are such that $f(f(m)) \le \min
I$ and $\# I \le m$, then
$$\left\|\sum_{i\in I} \lambda_i z_i\right\|_S \ge \frac{G(m)}{f(m)}
\sum_{i\in I} |\lambda_i|,
\text{\ whenever\ } (\lambda_i)_{i\in I}\subset \R.$$
\end{lem}

\begin{proof}
Using (\ref{eq2.10}) we can write
 for any $n$ and $\ell = 0,1,\ldots, L_n+1$ $x^*_n$ as

$$x^*_n= \sum_{\stackrel{\mu\in\mathcal{T}}{|\mu|=\ell}} \beta(\mu)x^*_{n,\mu}$$
where $\beta(\mu) = 1/\prod^{|\mu|-1}_{i=1} f(k_{(\mu_1,\ldots,\mu_i)})$
and $x^*_{n,\mu}$ is associated to 
$\{\nu \in\cup^{L_n+1-\ell}_{j=0} \N^j\mid (\mu,\nu)\in\mathcal{T}\}$, 
for $\mu\in\mathcal{T}$, with $|\mu|=\ell$, and where
$(x^*_{n,\mu})_{\mu\in \mathcal{T}, |\mu|=\ell}$ is a block with respect to
$\prec_{\text{lex}}$. Assume $m\in \N$, $I\subset \N$, with
$f(f(m)) \le \min I$ and $\# I \le m$, and let $(\lambda_i)_{i\in I} \subset
\R$. Without loss of generality we can assume that $\lambda_i \ge 0$,
for $i\in I$.

By recursion we will choose for each $\ell=0,1,\ldots, L_{\lfloor 
f(f(m))\rfloor}$ a node
$\nu_\ell\in \mathcal{T}$, with $|\nu_\ell|=\ell$ and
$\nu_0 \prec \nu_1 \prec \nu_2 \prec\ldots\prec \nu_\ell$, so that
\begin{equation}\label{eq2.20}
\sum_{i\in I} \lambda_i x^*_{i,\nu_\ell}(z_i) \ge \prod^{\ell-1}_{s=1}
\frac{f(k_{\nu_s})}{k_{\nu_s}} \sum_{i\in I} \lambda_i.
\end{equation} 
For $\ell=0, \mu_0 = \emptyset$, and since $x^*_{n,\emptyset} = x^*_n$, 
the claim
follows from the assumptions on $(z_n)$. Assume $\nu_0 \prec \nu_1 \prec
\ldots \prec \nu_\ell$, with $\ell < L_{\lfloor f(f(m))\rfloor}$, have been chosen.
Then it follows from (\ref{eq2.10}) that
$$x^*_{n,\nu_\ell} = \frac1{f(k_{\nu_\ell})} \sum^{k_{\nu_\ell}}_{j=1}
x^*_{i,(\nu_\ell,j)}$$
and by the induction hypothesis that
$$
\prod^{\ell-1}_{s=0} \frac{f(k_{\nu_s})}{k_{\nu_s}} \sum_{i\in I} \lambda_i
\le \sum_{i\in I} \lambda_i x^*_{i,\nu_\ell} (z_i)
= \frac1{f(k_{\nu_\ell})} \sum_{i\in I} \sum^{k_{\nu_\ell}}_{j=1} \lambda_i
x^*_{i,(\nu_\ell,j)}(z_i).
$$
Thus there exists a $j\le k_{\nu_\ell}$ so that for $\nu_{\ell+1} :=
(\nu_\ell,j)$ it follows that
$$\sum_{i\in I} \lambda_i x^*_{i,\nu_{\ell+1}} (z_i) \ge \prod^\ell_{s=0}
\frac{f(k_{\nu_s})}{k_{\nu_s}} \sum_{i\in I} \lambda_i.$$
This finishes the proof of the claim.

For $\ell=0,1,\ldots, L_{\lfloor f(f(m))\rfloor}$ let 
$$
y^*_\ell = \frac1{f(mk_{\nu_\ell})} \sum_{i\in I} \sum^{k_{\nu_\ell}}_{j=1}
x^*_{i,(\nu_\ell,j)}
= \frac{f(k_{\nu_\ell})}{f(mk_{\nu_\ell})} \sum_{i\in I} x^*_{i,\nu_\ell}
$$
(recall that $x^*_i$, with $i\in I$, is the associated vector to a tree of 
length at
least $L_{\lfloor f(f(m))\rfloor}+1$). Since $\# I \le m$ it follows that on the
$\|y^*_\ell\|_S \le 1$,
and, then, that
\begin{align*}
\left\|\sum_{i\in I} \lambda_iz_i\right\|_S &\ge
\max_{\ell\le L_{\lfloor f(f(m))\rfloor}} y^*_\ell \left(\sum_{i\in I}
\lambda_iz_i\right)
=\max_{\ell\le L_{\lfloor f(f(m))\rfloor}} \frac{f(k_{\nu_\ell})}{f(m k_{\nu_\ell})}
    \sum_{i\in I}\lambda_i x^*_{i,\nu_\ell}(z_i)\\
 &\ge \max_{\ell\le L_{\lfloor f(f(m))\rfloor}}\frac{f(k_{\nu_\ell})}{f(m k_{\nu_\ell})} \prod^{\ell-1}_{s=1}
\frac{f(k_{\nu_s})}{k_{\nu_s}} \sum_{i\in I} \lambda_i
 = \frac{G(m)}{f(m)} \sum_{i\in I} \lambda_i.
\end{align*}
In the last of above  inequalities we used (\ref{eq2.19}), the  remark after 
(\ref{eq2.19}),
 and Proposition~\ref{pro2.2}.
\end{proof} 

Before we can prove that our chosen sequence $(x^*_n)$ satisfies the
properties stated in Proposition~\ref{pro2.2} we will need another argument.
For $r\ge 2$ and $x\in S$ define
$$|||x|||_r = \sup_{r \leq \ell <\infty} \|x\|_\ell = 
\sup_{\ell\ge r, E_1<E_2<\cdots<
E_\ell} \frac1{f(\ell)} \sum^\ell_{i=1} \|E_ix\|_S.$$
 Note that $|||\cdot|||_r$ is an equivalent norm on $S$. The next lemma makes 
the
following qualitative statement precise:

If $x\in c_{00}$ and if $r\ge 2$ is ``big enough'', if $\ell\ge r$ and $E_1 <
E_2 <\cdots< E_\ell$ are such that
$$|||x|||_r = \|x\|_\ell = \frac1{f(\ell)} \sum^\ell_{i=1} \|E_ix\|_S$$
then ``for most of the'' $i\in \{1, \ldots ,\ell\}$ it follows that either
$\|E_ix\|_S = \|E_ix\|_{\ell_\infty}$ or $\|E_ix\|_S = \|E_ix\|_{n_i}$ with 
$n_i$ 
being
``much bigger than $r$''.

\begin{lem}\label{lem2.6}
There is a $d>1$, so that for all $r$ with $f(r) >d^2$ and all 
$x= \sum x_j e_j \in c_{00}$, it follows that
 $$|||x|||_r = \|x\|_\ell \le
\frac1{1-\frac{d}{\sqrt{f(r)}}} \frac1{f(\ell)}
\left(\sum_{j\in J} |x_j| + \sum^{\ell-\# J}_{i=1}
\|E_ix\|_{n_i}\right),$$
where $\ell\in\N$, $\ell\ge r$, $J \subset\text{\rm supp}(x)$, 
 $\# J \le \ell,\quad E_1 < E_2 < \cdots < E_{\ell-\# J}$,
$E_i\subset \N$ and $E_i\cap J=\emptyset$,  and  $n_1,n_2,\ldots, n_{\ell-\# 
J}\in\N$,
  $n_j > r^{f(r)}$ for $j=1,2,\ldots, \ell-\# J$.
\end{lem}

Lemma \ref{lem2.6} is almost identical to Lemma~5 in \cite{S2}; for the sake
of being self-contained we include a proof at the end of this section.
Now we are ready for the last step of proving Theorem~\ref{thm1.2}, by showing
that $(x^*_n)$ as chosen before satisfies Proposition~\ref{pro2.2}.

\begin{lem}\label{lem2.7}
There is a constant $c>1$ and an $m_0\in \N$, so that for any block
sequence $(z_n)\subset S$, with $x^*_n(z_n) = 1$ and 
$x^*_{i-1}<z_i<x^*_{i+1}$, for $n\in \N$, and any $m\ge m_0$ it
follows that
$$\left\|\sum^\infty_{i=1} \lambda_i e_i\right\|_m \le \frac{c}{\sqrt{G(m)}}
\left\|\sum^\infty_{i=1}\lambda_i z_i\right\|_S$$
whenever $(\lambda_i) \in c_{00}$. 
\end{lem}

Choosing now $C(\ell) = \sqrt{G(\ell)}/c$ for $\ell\ge \min\{m \ge
m_0\colon \sqrt{G(m)}\ge c\}$, and $C(\ell)=1$ for $\ell <
\min\{m\ge m_0 : \sqrt{G(m)}\ge c\}$ we note that the
assumptions of Proposition~\ref{pro2.2} are fulfilled. For the case 
$\ell<\min\{m\ge m_0 : \sqrt{G(m)}\ge c\}$ we are simply using the fact that 
every
 block basis in $S$ whose elements are of norm at least 1 dominates  $(e_i)$ [S2].

\begin{rem} It will be important for the arguments in Section \ref{sec3} that 
 for any $r>1$ the series $\sum 1/C(r_\ell)$
  converges, where $r_0=r$ and, inductively, $r_{\ell+1}=r_{\ell}^{f(r_{\ell})}$
  (Lemma~\ref{lem2.4}).
\end{rem}

\begin{proof}[Proof of Lemma~\ref{lem2.7}]
We choose $m_0$ so that $f(m_0)\ge 2d^2$, where $d$ is chosen as in
Lemma~\ref{lem2.6}, and so that $m_0>k_1$, the first element of the previously
chosen sequence $(k_n)$. For $r\in \R_+$ with $f(r)\ge
2d^2$ we consider the sequence $(r_\ell)^\infty_{\ell=0}$, with $r_0=r$ and
$r_{\ell+1}  = r^{f(r_\ell)}_\ell$ and observe that the two series
$\sum^\infty_{\ell=0} \frac{d}{\sqrt{f(r_\ell)}-d}$ and
$\sum^\infty_{\ell=0} \frac{f(f(r_\ell))\sqrt{G(r_\ell)}}{f(r_\ell)}$
are finite (note that $\sqrt{G(r_\ell)} \le \sqrt{f(r_\ell)}$ and that
$\frac{f(f(r_\ell))}{\sqrt{f(r_\ell)}}$ is summable).
By Lemma~\ref{lem2.4} also the series 
$\sum^\infty_{\ell=0}\frac1{\sqrt{G(r_\ell)}}$ 
is finite. Therefore
$$c(r) := \prod^\infty_{\ell=0} \frac1{1-\frac{d}{\sqrt{f(r_\ell)}}} \left[ 1
+  \frac{f(f(r_\ell))\sqrt{G(r_\ell)}}{f(r_\ell)} + 
\frac1{\sqrt{G(r_\ell)}-1}\right]$$
 is uniformly bounded on $[m_0,\infty)$. By  induction on $\# I$ we will show 
that
\begin{equation}\label{eq2.21}
\left\|\sum_{i\in I} \lambda_i e_i\right\|_m \le \frac{c(m)}{\sqrt{G(m)}}
\left\|\sum_{i\in I} \lambda_iz_i\right\|_S \text{ for all } m\ge m_0.
\end{equation}
This would imply the statement of the Lemma if we choose $c = \sup_{m\ge m_0} 
c(m)$.

 If $\# I=1$ the claim is trivial (note that $\|e_i\|_m=1/f(m) \ge 
1/\sqrt{G(m)}$). 

Assume that for some $k\in
\N$ and all $I\subset \N$, $\# I\le k$, (\ref{eq2.21}) is true
and assume $\# I = k+1$, $m\ge  m_0$, and $(\lambda_i)_{i\in I} \subset
\R$.
By Lemma~\ref{lem2.6} we can choose $\ell\ge m$, $J \subset I$,
$\# J \le\ell$, $E_1 < E_2 <\cdots< E_{\ell-\# J}$, $E_i\cap J=\emptyset$,
and $n_1,n_2,\ldots, n_{\ell-\# J} \in (\lfloor m^{f(m)}\rfloor ,\infty) \cap\N$ 
so that (put $x  = \sum_{i\in I} \lambda_i e_i$ and $z =\sum_{i\in I} 
\lambda_iz_i$)
\begin{align*}
\left\|\sum_{i\in I} \lambda_i e_i\right\|_m &\le \left|\left|\left|
\sum_{i\in I} \lambda_i e_i\right|\right|\right|_m
= \left\|\sum_{i\in I} \lambda_i e_i\right\|_\ell\\
&\le \frac1{1-\frac{d}{\sqrt{f(m)}}} \frac1{f(\ell)} \left[\sum_{j\in
J} |\lambda_j| + \sum^{\ell-\# J}_{i=1}
\|E_ix\|_{n_i}\right]\\
&\le \frac1{1-\frac{d}{\sqrt{f(m)}}} \frac1{f(\ell)} \left[\sum_{
\stackrel{j\in J}{j\le f(f(\ell))}} |\lambda_i| +
\sum_{\stackrel{j\in J}{j>f(f(\ell))}} |\lambda_i| + \sum^{\ell-
\# J}_{i=1} \|E_ix\|_{n_i}\right]\\
&\le \frac1{1-\frac{d}{\sqrt{f(m)}}} \left[\frac{f(f(\ell))}{f(\ell)}
\max_{\lambda\in I} |\lambda_i| + \frac1{f(\ell)} \left(\sum_{\stackrel{j\in
J}{j>f(f(\ell))}} |\lambda_i| + \sum^{\ell-\# J}_{i=1}
\|E_ix\|_{n_i}\right)\right]\\
&\le \frac1{1-\frac{d}{\sqrt{f(m)}}} \left[ \frac{f(f(m))}{f(m)} \|z\|_s +
\frac1{f(\ell)} \left(\sum_{\stackrel{j\in J}{j>f(f(\ell))}}
|\lambda_i| + \sum^{\ell-\# J}_{i=1} \|E_i x\|_{n_i} \right)\right].
\end{align*}

We can assume that $E_i\subset \supp(x)$ for $i\le\ell-\#J$.
If $J\not=\emptyset$ then the cardinality of the support of $E_i x$ is smaller 
than the 
cardinality of supp$(x)$.  If $J=\emptyset$  and if, say, 
$E_1=\text{supp}(x)$ then we could
  split $E_1$ into $\tilde E_1$ and $\tilde E_2$, choose $n_2=n_1$ and
observe that above inequalities
   still hold.  Thus we can assume  that for all $i=1,\ldots \ell-\#J$, the 
support of $E_i(x)$ is
    of lower cardinality than the support of $x$. Thus we can assume that our induction hypothesis
     applies to $E_i(x)$, $i\le\ell-\#J$.

We distinguish now between two cases.

If $\sum\limits_{\stackrel{j\in J}{j>f(f(\ell))}} |\lambda_j| \ge
\frac1{\sqrt{G(\ell)}-1} \sum\limits^{\ell-\# J}_{i=1}
\|E_i(x)\|_{n_i}$, then it follows from Lemma~\ref{lem2.5}  (applied to $\ell=m$ and
 $I=\{j\in J: j>f(f(\ell))\}$) that
\begin{align*}
\frac1{f(\ell)} \left(\sum_{\stackrel{j\in J}{j>f(f(\ell))}}
|\lambda_j| + \sum^{\ell-\# J}_{i=1} \|E_ix\|_{n_i}\right)
& \le \frac1{f(\ell)} \sqrt{G(\ell)} \sum_{\stackrel{j\in J}{j > f(f(\ell))}} 
|\lambda_j|\\
& \le \frac1{\sqrt{G(\ell)}} \left\|\sum_{\stackrel{j\in J}{j
>f(f(\ell))}} \lambda_jz_j\right\|_S \quad(\#J\le\ell)\\
& \le \frac1{\sqrt{G(m)}} \left\|\sum_{i\in I} \lambda_iz_i\right\|_S.
\end{align*}

If $\sum\limits_{\stackrel{j\in J}{j>f(f(\ell))}} |\lambda_j| <
\frac1{\sqrt{G(\ell)}-1} \sum\limits^{\ell-\# J}_{i=1}
\|E_i(x)\|_{n_i}$ we deduce that
\begin{align*}
\frac1{f(\ell)} \Biggl( \sum_{\stackrel{j\in J}{j\ge f(f(\ell))}}|\lambda_j|  
 &+  \sum^{\ell-\# J}_{i=1} \|E_i(x)\|_{n_i} \Biggr) \\
&\le \frac1{f(\ell)} \left(1 + \frac1{\sqrt{G(\ell)}-1}\right)
\sum^{\ell - \# J}_{i=1} \|E_i(x)\|_{n_i}\\
&\le \left(1 + \frac1{\sqrt{G(m)}-1}\right) c(m^{f(m)})
\frac1{\sqrt{G(m^{f(m)})}} \frac1{f(\ell)} \sum^{\ell-\# J}_{i=1}
\|F_i(z)\|_S \\
&\text{(by the induction hypothesis applied to $E_ix$)}\\
\intertext{here we let $z = \sum\limits_{i\in I} \lambda_iz_i$ and
$F_i = \bigcup\limits_{j\in E_i} \text{ supp}(z_j)$, and note
that $F_iz = \sum\limits_{j\in E_i} \lambda_jz_j$.}
&\le \left(1 + \frac1{\sqrt{G(m)}-1}\right) c(m^{f(m)})
\frac1{\sqrt{G(m)}} \|z\|_S.
\end{align*}
Thus in both cases we conclude that 
\begin{align*}
\left\|\sum \lambda_i e_i\right\|_m &\le
\left(\frac1{1-\frac{d}{f(m)}}\right) \left(\frac{f(f(m))\cdot
\sqrt{G(m)}}{f(m)} + 1 + \frac1{\sqrt{G(m)}-1}\right) c(m^{f(m)})
\frac1{\sqrt{G(m)}} \|z\|_s\\
&= c(m) \frac1{\sqrt{G(m)}} \left\|\sum_{i\in I}\lambda_i z_i\right\|_S,
\end{align*}
which finishes the proof of Lemma~\ref{lem2.7}.
\end{proof}

We still have to prove Lemma~\ref{lem2.4} and Lemma~\ref{lem2.6}.

\begin{proof}[Proof of Lemma~\ref{lem2.4}]
Note that for $r>1$
\begin{equation}\label{eq2.22}
G(r) \ge \widetilde G(r) := \sum^\infty_{j=1} 1_{[m_{j-1},m_j)}(r)
 \frac{f(r) f(k_j)}{f(r k_j)}\prod^{j-1}_{s=1} \frac{f(k_s)}{k_s} \qquad
(m_0=1).
\end{equation}
Therefore it will be enough to show that for $j\ge 2$,
\begin{equation}\label{eq2.23}
\sum_{  \ell\in \N , m_{j-1}\le r_\ell <m_j  }
\frac1{\sqrt{\tilde G(r_\ell)}} \le \vp_{j-1} + \vp_j.
\end{equation}
Using (\ref{eq2.12}) we can easily prove by induction on $\ell\in \N_0$
that
\begin{equation}\label{eq2.24}
r_\ell \ge r^{\left(\frac34 f(r)\right)^{2^\ell-1}} \ge m_{j-1}^{\left(\frac34 
f(m_{j-1})\right)^{2^\ell-1}}
\end{equation}
and (by applying $f$) it follows from (\ref{eq2.13}) that
 (recall that $m_{j-1}\ge m_1\ge k_1$)\begin{equation}\label{eq2.25}
f(r_\ell) \ge \left(\frac34 f(m_{j-1})\right)^{2^\ell}.
\end{equation}
We distinguish between two kinds of $r_\ell$'s in $\{r_\ell\colon \ m_{j-1}
\le r_\ell < m_j\}$. If $0\le \ell$ is such that $m_{j-1} \le r_\ell \le k_j$
we deduce that 
\begin{align}\label{eq2.26}
\widetilde G(r_\ell) &= \frac{f(r_\ell)f(k_j)}{f(r_\ell k_j)}\prod^{j-1}_{s=1} 
\frac{f(k_s)}{k_s} 
  \ge  f(r_\ell) \frac{f(k_j)}{f(k^2_j)}\prod^{j-1}_{s=1} \frac{f(k_s)}{k_s}\\
&\ge \frac12  f(r_\ell)\prod^{j-1}_{s=1} \frac{f(k_s)}{k_s} \ge 
\frac12 \left(\frac34 f(m_{j-1})\right)^{2^\ell}\quad\prod^{j-1}_{s=1} 
\frac{f(k_s)}{k_s} \quad\
                   \text{(by (\ref{eq2.25}))}\nonumber\\
&\ge \frac3{16} k_{j-1} \left(\frac38 k_{j-1}\right)^{2^\ell-1}\quad 
\prod^{j-1}_{s=1} \frac{f(k_s)}{k_s}\nonumber\\
\intertext{(note that $f(m_{j-1})\ge f(m_{j-1} k_{j-1}) - f(k_{j-1}) \ge
\frac12 k_{j-1}$)}
&= \frac3{16} f(k_{j-1}) \left(\frac38 k_{j-1}\right)^{2^\ell-1}\quad
          \prod^{j-2}_{s=1} \frac{f(k_s)}{k_s} 
 \ge \frac3{16}f(k_{j-1}) \left(\frac38 k_1\right)^\ell 
   \prod^{j-2}_{s=1} \frac{f(k_s)}{k_s} .\nonumber
\end{align}
Thus
\begin{equation}\label{eq2.27}
\sum_{\stackrel{\ell\in \N_j}{m_{j-1}\le r_\ell\le k_j}} 
\frac1{\sqrt{\widetilde{G}(r_\ell)}} 
\le \left(\frac{16}3 
\frac1{f(k_{j-1})}\prod^{j-2}_{s=1} \frac{k_s}{f(k_s)}\right)^{1/2} 
\sum^\infty_{\ell=0} \left(\frac83
k_1\right)^{-\ell/2}
< \vp_{j-1}\qquad \text{(by (\ref{eq2.14})).}
\end{equation}
To estimate $1/\sqrt{\widetilde G(r_\ell)}$, with $k_j < r_\ell < m_j$, 
we put
 $\ell_0 = \min\{\ell\mid r_{\ell_0}>k_j\}$ and $\ell_1 =
\max\{\ell\mid r_\ell < m_j\}$
and observe for $\ell\in \{\ell_0, \ell_0+1,\ldots, \ell_1\}$ that
\begin{equation}\label{eq2.28}
\widetilde G(r_\ell) =  \frac{f(r_\ell)f(k_j)}{f(r_\ell
 k_j)}\prod^{j-1}_{s=1} \frac{f(k_s)}{k_s}
\!\ge\! f(k_j) \frac{f(r_\ell)}{f(r^2_\ell)}\prod^{j-1}_{s=1} \frac{f(k_s)}{k_s} 
\!\ge\! \frac12  f(k_j)\prod^{j-1}_{s=1} \frac{f(k_s)}{k_s}.
\end{equation}
As in (\ref{eq2.24}) and (\ref{eq2.25}) we observe that
 $$r^{\left(\frac34 f(r_{\ell_0})\right)^{2^{\ell_1-\ell_0}-1}}_{\ell_0} \le
r_{\ell_1} < m_j
\text{\ and\ }
  \left(\frac34 f(r_{\ell_0})\right)^{2^{\ell_1-\ell_0}} \le f(m_j) \le k_j$$
which implies by (\ref{eq2.13}) that
 $2^{\ell_1-\ell_0} f(\frac34 f(r_{\ell_0})) \le \frac43 f(k_j)$
and thus (by (\ref{eq2.12}))
\begin{equation}\label{eq2.29}
\ell_1-\ell_0 \le \log_2 \frac{f(k_j)}{\frac34f\left(\frac34
f(r_{\ell_0})\right)} \le \log_2 f(k_j).
\end{equation}
The inequalities (\ref{eq2.28}), (\ref{eq2.29}) and (\ref{eq2.14}) imply that
\begin{align}\label{eq2.30}
\sum_{\ell\in \N, k_j < r_\ell < m_j} \widetilde G(r_\ell)^{-1/2} &\le
2(\ell_1-\ell_0) \left( \frac1{f(k_j)}\prod^{j-1}_{s=1} \frac{k_s}{f(k_s)}
\right)^{1/2}\\
&\le 2\frac{\log_2 f(k_j)}{\sqrt{f(k_j)}} \left(\prod^{j-1}_{s=1}
\frac{k_s}{f(k_s)}\right)^{1/2} \le \vp_j.\nonumber
\end{align}

Finally (\ref{eq2.27}) and (\ref{eq2.30}) imply the claimed inequality
(\ref{eq2.23}).
\end{proof}

\begin{proof}[Proof of Lemma \ref{lem2.6}]
Since $f$ is close to a logarithmic function we can choose a $c>1$ so that
for $\xi,\xi', R\ge 1$ and $0 < r < 1$
\begin{align}
\label{eq2.31}
\frac1c (f(\xi) + f(\xi')) &\le f(\xi\cdot\xi') \le f(\xi) + f(\xi')\\
\label{eq2.32}
\frac1c Rf(\xi) &\le f(\xi^R) \le Rf(\xi)\\
\label{eq2.33}
rf(\xi) &\le f(\xi^r)\\
\label{eq2.34}
\frac1c f(\xi) &\le f(\xi) -1\quad \text{if}\quad \xi\ge 2.
\end{align}
Assuming furthermore that $c\ge 2$ we deduce  for $\xi\ge 1$ that
\begin{equation}
\label{eq2.35}
f(\xi^{1/\sqrt{f(\xi)}}) = \log_2 (\xi^{1/\sqrt{f(\xi)}} +1)
\le 1 + \frac1{\sqrt{f(\xi)}} \log_2(\xi)
\le 1 + \sqrt{f(\xi)} \le c\sqrt{f(\xi)}, .
\nonumber
\end{equation}
Let $d = 4c^3$, let $r\in \R$, such that $f(r)>d^2$, and $x\in
c_{00}$. Choose $r\le \ell < \infty$ and $E_1 < E_2 <\cdots< E_\ell$ so that
  $|||x|||_r = \|x\|_\ell = \frac1{f(\ell)} \sum^\ell_{i=1} \|E_i(x)\|_S$. 
For $i\in \{1,2,\ldots, \ell\}$ let $n_i\in \N\cup \{\infty\}$ so that
$||E_ix\|_S = \|E_ix\|_{n_i}$. For two numbers $\tilde r, \widetilde R$,
$2\le \tilde r < \widetilde R < \infty$, we let
 $M := M(\tilde r,\widetilde R) := \{i\in \ell\colon \ \tilde r \le n_i <
\widetilde R\}$
and we choose for $i\in M$, $E^{(i)}_1 <\cdots< E^{(i)}_{n_i}$, with
$E^{(i)}_j \subset E_i$ whenever $j\le n_i$ so that
 $\|E_ix\|_S = \frac1{f(n_i)} \sum^{n_i}_{j=1} \|E^{(i)}_jx\|_S$. 
Now we observe that 
$\{E_i\colon \ i\notin M\} \cup \bigcup\nolimits_{i\in M} \{E^{(i)}_j\colon 1\le 
j\le n_i\}$
is well ordered by $<$ and its cardinality is $\ell - \# M + \sum_{i\in M}n_i$ 
which is at least $\ell$ and at most $\widetilde R \ell$. Thus we deduce
\begin{align*}
|||x|||_r &= \frac1{f(\ell)} \sum^\ell_{i=1} \|E_ix\|_S
\ge \frac1{f\left(\ell-\# M+ \sum\limits_{i\in M}n_i\right)}
\left[\sum^\ell_{\stackrel{i=1}{i\notin M}} \|E_ix\|_s + \sum_{i\in M}
\sum^{n_i}_{j=1} \|E^{(i)}_jx\|_S\right]\\
&\ge \frac1{f(\ell\cdot\widetilde R)} \left[\sum^\ell_{i=1, i\notin M}
\|E_ix\|_S + \sum_{i\in M} f(n_i) \|E_ix\|_S\right]
\mbox{ (since $n_i \le \widetilde R$ for $i\in M$)}\\
&\ge \frac1{f(\ell\cdot\widetilde R)} \left[\sum^\ell_{i=1} \|E_ix\|_S +
\sum_{i\in M} (f(\tilde r)-1) \|E_ix\|_S\right]\\
&\ge \frac1{f(\ell\cdot\widetilde R)} \left[\sum^\ell_{i=1} \|E_ix\|_S +
\frac{f(\tilde r)}c \sum_{i\in M} \|E_ix\|_S\right].
\end{align*}
Solving these inequalities for $\frac1{f(\ell)} \sum\limits_{i\in M}
\|E_ix\|_S$ we obtain that
\begin{align*}
\frac1{f(\ell)} \sum_{i\in M} \|E_ix\|_S &\le \frac1{f(\ell)}
\left[\frac1{f(\ell)} - \frac1{f(\ell\cdot\widetilde R)} \right]
\frac{cf(\ell\widetilde R)}{f(\tilde r)} \sum^\ell_{i=1} ||E_ix\|_S\\
&= \frac{c}{f(\tilde r)} \left[\frac{f(\ell\cdot\widetilde R) -
f(\ell)}{f(\ell)}\right] |||x|||_r\\
&\le c \frac{f(\widetilde R)}{f(\tilde r)  f(\ell)} |||x|||_r
  \le c \frac{f(\widetilde R)}{f(\tilde r) f(r)}|||x|||_r.
\end{align*}
Choosing for the $(\tilde r, \widetilde R)$ the values $(2,
r^{1/\sqrt{f(r)}})$, $(r^{1/\sqrt{f(r)}},r)$, $(r,r^{\sqrt{f(r)}})$ and
$(r^{\sqrt{f(r)}}, r^{f(r)})$ we deduce for all choices
 $\frac{f(\widetilde R)}{f(\tilde r) f(r)} \le \frac{c^2}{\sqrt{f(r)}}$.
 This implies that
  $$\frac1{f(\ell)} \sum_{\stackrel{1\le i\le \ell}{n_i<r^{f(r)}}} \|E_i x\| \le
    \frac{4c^3}{\sqrt{f(r)}} |||x|||_r$$
 and thus that
$$
|||x|||_r = \frac1{f(\ell)} \left[\sum^\ell_{\stackrel{i=1}{n_i<r^{f(r)}}}
\|E_ix\|_{n_i}  + \sum^\ell_{\stackrel{i=1}{n_i>r^{f(r)}}} \|E_ix\|_{n_i}\right]\\
 \le \frac{d}{\sqrt{f(r)}} |||x|||_r + \frac1{f(\ell)}
\sum^\ell_{\stackrel{i=1}{n_i>r^{f(r)}}} \|E_ix\|.
$$

 Now let $I = \{i\le \ell\mid n_i=\infty \text{ or } n_i\le r^{f(r)}\}$ and choose for $i\in I$ a
 $j_i\in \supp(x)\cap E_i$ so that $\|E_ix\|_{\ell_\infty} = |x_{j_i}|$, and
reorder
$(n_i)_{i\in \{1\ldots\ell\} \backslash I}$ and $(E_i)_{i\in \{1\ldots \ell\}
\backslash I}$ into $\tilde n_1,\tilde n_2,\ldots, \tilde n_{\ell- \# I}$
and $\widetilde E_1 < \widetilde E_2 <\cdots< \widetilde
 E_{\ell - \# I}$. Letting $J = \{j_i\colon \ i\in I\}$ we obtain the claimed
inequality:
$$|||x|||_r \le \frac1{1-\frac{d}{\sqrt{f(r)}}}  \frac1{f(\ell)}
\left[\sum_{j\in J} |x_j| + \sum^{\ell- \# J}_{i=1}
\|\widetilde E_ix\|_{\tilde n_i}\right].$$

\end{proof}
\section{The construction of strictly singular non-compact operators
on $GM$}\label{sec3}

The goal of this section is to prove Theorem \ref{thm1.1}.
 We postpone the
definition of the space $GM$ and only state some properties needed for the proof
 of Theorem \ref{thm1.1}.

 Recall that the {\em spreading model} of a seminormalized weakly null
sequence $(x_i)$ in some Banach space $(X, \| \cdot \|)$ is a sequence
$(y_i)$ along with a norm $| \cdot |$ on the span of $(y_i)$ such that
$$
\left| \sum_{i=1}^k \lambda_i y_i \right| =
\lim _{\substack{N \to \infty \\ N \leq n_1<n_2< \cdots n_k}}
\| \sum_{i=1}^k \lambda_i x_i \|
$$
for every finite sequence of scalars $(\lambda_i)_{i=1}^k$.

\begin{pro} \label{T:spreading}
The spreading model of the unit vector basis of $GM$ is the unit
vector basis of $S$.
\end{pro}

 In order to
defining $GM$ a lacunary set $J \subseteq \N$ is used which has the property  that if $n,m
\in J$ and $n<m$ then $4n^2 \leq \log \log \log m$, and $f(\min J)
\geq 256$.

\begin{pro} \label{T:GM}
 The norm of $GM$ satisfies
\begin{equation} \label{ena}
 \|x\|_S\le \| x  \|_{GM} \leq \| x \|_{S} + \sum_{\ell \in J} \| x \|_\ell, \text{ for }x\in c_{00}
\end{equation}
(here $\|\cdot\|_\ell$ still denotes the equivalent norm on $S$ introduced at the end of Section
 \ref{sec1}).
\end{pro}
Propositions \ref{T:spreading} and \ref{T:GM} will be shown at end of this section.
 The first inequality in (\ref{ena}) is an immediate consequence of the definition of $GM$ and
 the second inequality  will follow from Lemma
 \ref{L:SGM}.

Recall from the previous section that
there exists a seminormalized
block  sequence $(x_i^*)$ in $S^*$ and
a non-decreasing function $C:[2,\infty) \to (0,\infty)$
which  satisfies (\ref{eq2.1}).
Secondly we observe that from the remark following Lemma \ref{lem2.7}
and the conditions on $J$ it follows that

\begin{equation} \label{E:sumCJ}
\sum_{\ell \in J} \frac{1}{C(\ell)} <\infty
.\end{equation}

The sequence $(x_i^*)$ may not be seminormalized in
$GM^*$, i.e. $(x_i^*)$ could be a null sequence in
$GM$ (which would imply that the operator $T$
defined in Section \ref{sec2} would be compact on $GM$).
But Proposition  \ref{T:spreading} implies that
we can replace each $x^*_i$ by $\tilde x^*_i$ which has
the same distribution and up to some arbitrarily small number $\vp>0$
the same norm in $GM^*$ as in $S^*$. Indeed,
if $z^*=\sum_{i=1}^k \lambda_i e^*_i \in S^*\cap c_{00}$ and $z=\sum_{i=1}^k  \mu_i e_i \in S$,  with $\|z\|_S=1$
and $z^*(z)=\|z\|_{S^*}$ it follows that
\begin{align*}
\|z^*\|_{S}&=\sum_{i=1}^k \lambda_i \mu_i
          \le \liminf_{\substack{N \to \infty \\ N \leq n_1<n_2< \cdots n_k}}
 \left\|\sum_{i=1}^k \lambda_i  e_{n_i} \right\|_{GM}
            \cdot \left\|\sum_{i=1}^k \mu_i  e^*_{n_i} \right\|_{GM^*}  \\
         &= \liminf_{\substack{N \to \infty \\ N \leq n_1<n_2< \cdots n_k}}
  \left\|\sum_{i=1}^k \mu_i  e^*_{n_i} \right\|_{GM^*}.
\end{align*}

Thus we will assume from now on that   $(x_i^*)$  is also seminormalized in $GM^*$.

Then define the operator $T:GM \to GM$ by $T= \sum_i
 {x}_i^* \otimes e_i$. We now present
the 

\bigskip
\noindent
{\em Proof of Theorem \ref{thm1.1}} We first show that $T$ is
bounded. Every $x \in GM$ can be written as $x= \sum \lambda_i z_i$
with $x_{i-1}^*< z_i<x_{i+1}^*$, and $x_i^*(z_i)=1$ for all $i$. Thus
$Tx= \sum \lambda_i e_i$. By (\ref{ena}) we have that for any $x \in
c_{00}$, 
\begin{eqnarray} \label{E:prenormofTnu}
\| Tx \|_{GM} &\leq & \| Tx \|_S + \sum_{\ell \in J} \| \sum \lambda_i
e_i \|_{\ell} \\
& \leq & \| Tx \|_S + \sum_{\ell \in J} \frac{1}{C(\ell)} \| \sum
\lambda_i z_i \|_S \mbox{ (by property (\ref{eq2.1}))} \nonumber \\
& \leq & \left( \| T: S \to S \| + \sum_{\ell \in J} \frac{1}{C(\ell)}
\right) \| x \|_S . \nonumber
 \\& \leq & \left( \| T: S \to S \| + \sum_{\ell \in J} \frac{1}{C(\ell)}
\right) \| x \|_{GM} . \nonumber
\end{eqnarray}
Now we show that $T$ is strictly singular.
 From the definition of $(x^*_i)$ it follows that $T$ has an infinite dimensional
  kernel. Since  $T$ can be written as $T=\lambda+\tilde T$ with $\tilde T$ being
   strictly singular it follows that $\lambda=0$.

Finally, $T$ is not a compact operator since $({x}_i^*)$ is
seminormalized in $GM^*$. \hfill $\square$

Let $X$ be either the space $S$ or $GM$. 
For any sequence $\nu = (\nu_i) \in \ell_\infty$ define the
operator 
$$
T_\nu = \sum_i \nu_i x_i^* \otimes e_i :X \to GM. 
$$
The above proof shows that $T_{\nu}$ is a bounded strictly singular
operator  with
$$
\| \nu \|_\infty \inf_i \| x^*_i \|_X \leq 
\| T_\nu : X \to GM \| \leq 
\left( \| T:S \to S \| + \sum_{\ell \in J} 
\frac{1}{C(\ell)} \right) \| \nu \|_\infty .
$$
Therefore $\ell_\infty$ embeds in the space of operators from $X$ to
$GM$.  If $\nu \in \ell_\infty \backslash c_0$ then $T_\nu :X \to GM$ is 
non-compact. 

Let us now recall 
the definition of the space $GM$.

 Let ${\bf Q}$ be the set of scalar sequences
with finite support and rational coordinates  whose absolute value is at most one. 
 Write
 $J$ (introduced above) in increasing order as $\{ j_1, j_2, \ldots \}$. 
Now let $K \subset J$ be the set $\{ j_1, j_3, j_5 ,
\cdots \}$ and $L \subset J$ be the set $\{ j_2 , j_4 , j_6 , \ldots
\}$. Let $\sigma$ be an injective function from the collection of all
finite sequences of  elements of ${\bf Q}$ to $L$ such that if
$z_1, \ldots , z_i$ is such a sequence, then $(1/20)f(\sigma(z_1,
\ldots , z_i)^{1/40}) \geq \# \supp (\sum_{j=1}^i z_j)$. 
Then, recursively, we define a set
of functionals of the unit ball of $GM^*$ as follows: Let 
$$
GM^*_0= \{ \lambda e_n^* : n \in \N, \ | \lambda | \leq 1 \} . 
$$
Assume that $GM_k^*$ has been defined. Then $GM_{k+1}^*$ is the set of
all functionals of the form $E \ z^*$ where $E \subseteq \N$ is an interval
and $z^*$ has one of the following three forms:
\begin{equation}
z^*= \sum_{i=1}^\ell \alpha_i z_i^* 
\end{equation}
\hskip .4in where $\sum_{i=1}^\ell | \alpha_i | \leq 1$ and $z_i^* \in
GM_k^*$  for
$i=1, \ldots , \ell$.
\begin{equation}
z^*= \frac{1}{f(\ell)} \sum_{i=1}^\ell z_i^* 
\end{equation}
\hskip .4in where $z_i^* \in GM_k^*$ for $i=1, \ldots , \ell$, and
$z_1^*  < \cdots <
z^*_\ell$. 
\begin{equation}
z^*= \frac{1}{\sqrt{f(\ell)}} \sum_{i=1}^\ell z^*_i \mbox{ and }
z^*_i= \frac{1}{f(m_i)} \sum_{j=1}^{m_i} z^*_{i,j}
\end{equation}
\hskip .4in where $z^*_{i,j} \in GM^*_k$ for $1 \leq i \leq \ell$ and $1 \leq j
\leq m_i$, $z^*_{1,1} < \cdots < z^*_{1,m_1}< z^*_{2,1} < \cdots <
z^*_{\ell , m_\ell}$, 

\noindent
\hskip .4in $m_1= j_{2 \ell}$, $z_i^*$ has rational
coordinates, and $m_{i+1}= \sigma (z_1^*,
\ldots , z_i^*)$, for $i=1, \ldots , \ell -1$. 

Finally, the norm of $GM$ is defined by
$$
\| x \|_{GM} = \sup \{ z^* (x) : z^* \in \cup_{k=0}^\infty GM^*_k \}.
$$
For an interval $I \subseteq \N$ we define 
$$
J(I) = \{ \sigma(x_1^*, x_2^*, \ldots , x_n^*): n \in \N, \ 
x_1^* < x_2^* < \cdots < x_n^* , \ \min I \leq \max \supp (x_n^*) < \max I 
\}.
$$
Also, for $x^* \in GM^*$ we define $J(x^*)= 
J([ \min \supp (x^*), \max \supp (x^*) ] ).$

The next result relates
the  functionals of the unit ball of $GM^*$ to
the functionals of the unit ball of $S^*$.

\begin{lem} \label{L:SGM}
For any $z^* \in \cup_{k=0}^\infty GM_k$ there exist 
$T_0(z^*) \in \rm{Ba}\, (S^*)$, the unit ball of $S^*$, and a family 
$(T_j(z^*))_{j \in
J(z^*)} \subset {\rm Ba} \, S^*$ 
such that
\begin{itemize}
\item[1)] For $j \in \{ 0 \} \cup J(z^*)$,
$\supp T_j(z^*) \subseteq [ \min \supp z^* , \max \supp z^* ].$
\item[2)] For $j \in J(z^*)$
$$
T_j(z^*) \in \aco \bigl\{ \frac{1}{f(j)}: \sum_{s=1}^j x_s^* : x_1^* <
\cdots < x_j^* \mbox{ are in } {\rm Ba}\, (S^*) \bigr\}.
$$
where ``$\aco$'' denotes the absolute convex hull.
\item[3)] 
$$
z^*= T_0(z^*) +\sum_{j \in J(z^*)} T_j(z^*).
$$ 
\end{itemize}
\end{lem}

\begin{proof}
We proceed by induction on $k$ (assume that $z^* \in
GM^*_k$). For $k=0$, $z^*= \lambda e_n^*$ for some $\lambda \in \R$,
$| \lambda| \leq 1$ and some $n \in \N$. Then $J(z^*)= \emptyset$ and
$T_0(z^*)= z^*$. The inductive step, from $k$ to $k+1$ proceeds as
follows: By the definition of $GM_{k+1}^*$, we separate three cases:\\

\bigskip
\n{\em Case 1:} Assume that $z^*= E( \sum_{i=1}^\ell \alpha_i z_i^*)$ where
$E \subset \N$ is an interval, $z_i^* \in GM_k^*$ for all $i\le \ell$ and
$\sum_{i=1}^{\ell} | \alpha_i | \leq 1$. Let 
$\tilde{E} = [ \min \supp  E(z^*) , \max \supp E(z^*) ].$
Then $\tilde{E} \subseteq E$ and by the induction hypothesis we have
\begin{equation*}
z^*  =  \tilde{E} \left( \sum_{i=1}^\ell \alpha_i z_i^* \right)=
\sum_{i=1}^\ell \alpha_i \tilde{E} (z_i^* )
  =  \sum_{i=1}^\ell \alpha_i T_0( \tilde{E}(z_i^*)) + 
\sum_{i=1}^\ell \sum_{j \in J(\tilde{E} (z_i^*))} T_j(\tilde{E}(z_i^*)).
\end{equation*}
Set 
$$
T_0(z^*)= \sum_{i=1}^\ell \alpha_i T_0( \tilde{E}(z_i^*)) 
$$
and after noting that $\cup_{i=1}^\ell J(\tilde{E}(z_i^*)) \subseteq
J(z^*)$, for $j \in J(z^*)$ set
$$
T_j(z^*)= \sum_{ \stackrel{i=1, \ldots , \ell}{j \in J(\tilde{E}(z_i^*) }}
 \alpha_i T_j (\tilde{E} (z_i^*))
$$
where the sum over an empty set of indices is zero, and $T_j(z^*)=0$
if $j \in J(z^*) \backslash \cup_{i=1}^\ell J( \tilde{E}(z_i^*))$.
It is easy to see that the above choices of $T_0(z^*)$ and $T_j(z^*)$,
$j \in J(z^*)$, satisfy the conclusion of the lemma.

\bigskip
\n{\em Case 2:} Assume that $z^*= E \left( \frac{1}{f(\ell)}
\sum_{i=1}^\ell z_i^* \right)$ where $E \subseteq \N$ is an interval
and $z_1^*< z_2^* < \cdots < z_\ell^*$ in $GM_k^*$.  
By the induction hypothesis we have that
\begin{equation*}
z^* = \frac{1}{f(\ell)} \sum_{i=1}^\ell E(z_i^*) 
    =  \frac{1}{f(\ell)} \sum_{i=1}^\ell T_0(Ez_i^*) + 
\sum_{i=1}^\ell \sum_{j \in J(E(z_i^*))} \frac{1}{f(\ell)}T_j(Ez_i^*).
\end{equation*}
Set
$$
T_0(z^*)= \frac{1}{f(\ell)} \sum_{i=1}^\ell T_0(Ez_i^*)
$$
and after noting that $J(E(z_i^*)) \cap
 J(E(z_j^*))=\emptyset$ (since $\sigma$ is injective) for $1 \leq i \not = j \leq \ell$, and
$\cup_{i=1}^\ell J(E(z_i^*)) \subseteq J(z^*)$, for $j \in J(z^*)$ set
$$
T_j(z^*)= \frac{1}{f(\ell)}T_j(Ez_i^*) \mbox{ if }j \in J(E(z_i^*)),
$$
and $T_j(z^*)= 0$ if $j \in J(z^*) \backslash  \cup_{i=1}^\ell J(E (z_i^*))
$. It is easy to see that the conclusion of the
lemma is satisfied.

\bigskip
\n{\em Case 3:} Assume that $z^*= E \left( 
\frac{1}{\sqrt{f(\ell)}} \sum_{i=1}^\ell z_i^* \right)$ and $z_i^*=
\frac{1}{f(m_i)} \sum_{j=1}^{m_i} z_{i,j}^*$ 
where $z^*_{i,j} \in GM^*_k$ for $1 \leq i \leq \ell$ and $1 \leq j
\leq m_i$, $z^*_{1,1} < \cdots < z^*_{1,m_1}< z^*_{2,1} < \cdots <
z^*_{\ell , m_\ell}$, $m_1= j_{2 \ell}$, $z_i^*$ has rational
coordinates, and $m_{i+1}= \sigma (z_1^*,
\ldots , z_i^*)$, for $i=1, \ldots , \ell -1$.
Let 
$$
i_1= \min \bigl\{ i \in \{ 1, \ldots , \ell \} : E \cap \supp (z_i^*) \not =
\emptyset \bigr\}
$$
and 
$$
i_2 = \max \bigl\{ i \in \{ 1, \ldots , \ell \} : E \cap \supp (z_i^*) \not =
\emptyset \bigr\}.
$$
If $i_1=i_2$ then we proceed as in Case 2. Therefore without loss of
generality, we assume that $i_1<i_2$. 
Let 
$$
j_1= \min \bigl\{ j \in \{ 1, \ldots , m_{i_1} \}: E \cap \supp (z_{i_1,j}^*) \not =
\emptyset \bigr\}
$$
and
$$
j_2 = \max \bigl\{ j \in \{ 1, \ldots , m_{i_2} \bigr\}: E \cap \supp
(z_{i_2,j}^*) 
\not = \emptyset \} 
$$
By the induction hypothesis we have
\begin{eqnarray*}
E(z^*) &= & \frac{1}{\sqrt{f(\ell)}} \left[ 
\frac{1}{f(m_{i_1})} \sum_{j=j_1}^{m_{i_1}} E(z^*_{i_1,j}) 
+ \sum_{i=i_1+1}^{i_2-1} \frac{1}{f(m_i)} \sum_{j=1}^{m_i}z^*_{i,j}
+ \frac{1}{f(m_{i_2})} \sum_{j=1}^{j_2}E(z^*_{i_2,j}) \right] \\
& = & \frac{1}{\sqrt{f(\ell)}} \frac{1}{f(m_{i_1})} 
\sum_{j=j_1}^{m_{i_1}} T_0 (E(z_{i_1,j}^*)) \\
& & + \sum_{i=i_1+1}^{i_2-1} \frac{1}{\sqrt{f(\ell)}}\frac{1}{f(m_i)} 
\sum_{j=1}^{m_i} T_0(E z^*_{i,j}) 
+ \frac{1}{\sqrt{f(\ell)}}\frac{1}{f(m_{i_2})} 
\sum_{j=1}^{j_2} T_0(Ez^*_{i_2,j}) \\
& & + 
\sum_{j=j_1}^{m_{i_1}} \sum_{k \in J(Ez^*_{i_1,j})}
\frac{1}{\sqrt{f(\ell)}} \frac{1}{f(m_{i_1})} T_k(Ez^*_{i_1,j}) \\
& & + \sum_{i=i_1+1}^{i_2-1} \sum_{j=1}^{m_i} 
\sum_{k \in J(Ez^*_{i,j})} \frac{1}{\sqrt{f(\ell)}} \frac{1}{f(m_i)}  
T_k(E z^*_{i,j}) \\
& & + 
\sum_{j=1}^{j_2} \sum_{k \in J(Ez^*_{i_2,j})}
\frac{1}{\sqrt{f(\ell)}} \frac{1}{f(m_{i_2})} T_k(Ez^*_{i_2,j}) 
\end{eqnarray*}
Set
$$
T_0(Ez^*)= \frac{1}{\sqrt{f(\ell)}} \frac{1}{f(m_{i_1})} 
\sum_{j=j_1}^{m_{i_1}} T_0 (E(z_{i_1,j}^*)) 
$$
and after noting that 
\begin{equation} \label{E:union}
\{ m_{i_1+1}, \ldots , m_{i_2} \} \cup
\cup_{j=j_1}^{m_{i_1}} J(Ez^*_{i_1,j}) \cup 
\cup_{i=i_1+1}^{i_2-1} \cup_{j=1}^{m_i} J(E z^*_{i,j}) 
\cup \cup_{j=1}^{j_2} J(Ez^*_{i_2,j}) 
\subseteq J(Ez^*)
\end{equation}
and that the  sets $\{ m_{i_1+1}, \ldots , m_{i_2} \}$,
$J(Ez^*_{i_1,j})$ (for $j=j_1, \ldots ,m_{i_1}$), $J(E z^*_{i,j})$
(for $i=i_1+1, \ldots ,i_2-1$ and $j=1, \ldots ,m_i $),
$J(Ez^*_{i_2,j})$ (for $j=1, \ldots , j_2$)   
are mutually disjoint (by the injectivity of $\sigma$), set
$$
T_k(Ez^*) = \left\{ 
\begin{array}{ll}
\frac{1}{\sqrt{f(\ell)}}\frac{1}{f(m_i)} \sum_{j=1}^{k} T_0(Ez^*_{i,j})
& \mbox{ if } k =m_i \in \{ m_{i_1+1}, \ldots , m_{i_2 -1} \}\\
\frac{1}{\sqrt{f(\ell)}}\frac{1}{f(m_{i_2})}\sum_{j=1}^{j_2} T_0(Ez^*_{i_2,j})
& \mbox{ if }k =m_{i_2} \\
\frac{1}{\sqrt{f(\ell)}} \frac{1}{f(m_{i_1})} T_k(Ez^*_{i_1,j}) 
& \mbox{ if } k \in \cup_{i=j_1}^{m_{i_1}} J(Ez^*_{i_1,j}) \\
\frac{1}{\sqrt{f(\ell)}} \frac{1}{f(m_i)}  
T_k( Ez^*_{i,j}) & \mbox{ if }k \in \cup_{i=i_1+1}^{i_2-1} 
\cup_{j=1}^{m_i} J(z^*_{i,j})\\
\frac{1}{\sqrt{f(\ell)}} \frac{1}{f(m_{i_2})} T_k(Ez^*_{i_2,j}) 
& \mbox{ if } k \in \cup_{j=1}^{j_2}J(Ez^*_{i_2,j})\\
0 & \mbox{ if } k \in J(Ez^*) \mbox{ otherwise.}
\end{array} \right.
$$
It is easy to see that the conclusion of the lemma is
satisfied. \end{proof}

Note that inequality (\ref{ena}) (which was used in the proof of
Theorem \ref{thm1.1}) is an immediate consequence of Lemma
\ref{L:SGM}. It only remains to give the 

\bigskip
\noindent
{\em Proof of Proposition \ref{T:spreading}}.
Let $x = \sum_{i=1}^k \lambda_i e_i \in c_{00}$. We want to show that
\begin{equation} \label{E:spreading1}
\lim_{\substack{N \to \infty \\ N \leq n_1<n_2< \cdots <n_k}}
\| \sum_{i=1}^k \lambda_i e_{n_i}  \|_{GM}
= \| \sum_{i=1}^k \lambda_i e_i \|_S
.\end{equation}
Let $\e>0$. Since $J$ is lacunary enough and 
$\sigma$ is injective we can
choose  $N \in \N$
sufficiently large, such that
\begin{equation} \label{E:n1}
k \max_i | \lambda_i | \sum_{\ell \in J([N,\infty))} \frac{1}{f(\ell)} <
\e.
\end{equation}
Thus, if $N \leq n_1< n_2 < \cdots n_k$, then by Lemma \ref{L:SGM},
\begin{align*}
\| \sum_{i=1}^k \lambda_i e_{n_i} \|_{GM} &\! \leq \!  
\| \sum_{i=1}^k \lambda_i e_i \|_S + \sum_{ \ell \in J([N, \infty))}
\| \sum_{i=1}^k \lambda_i e_{n_i} \|_{\ell} \\
& \! \leq \!  \| \sum_{i=1}^k \lambda_i e_i \|_S + \sum_{i=1}^k | \lambda_i
| \sum_{\ell \in J([N, \infty))} \| e_{n_i} \|_{\ell} \\
& \! \leq \!  \| \sum_{i=1}^k \lambda_i e_i \|_S + k \max_i | \lambda_i |
\sum_{\ell \in J([N,\infty))}\!  \frac{1}{f(\ell)}
  \! \leq \!  \| \sum_{i=1}^k \lambda_i e_i \|_S +\e \mbox{ (by
(\ref{E:n1}))}
\end{align*}
which finishes the proof of (\ref{E:spreading1}). 
\hfill $\square$

{\footnotesize
\noindent
Department of Mathematics, University of South Carolina, Columbia, SC 29208, 
giorgis@math.sc.edu

\noindent
Department of Mathematics, Texas A{\&}M University, 
College Station TX 77843, schlump@math.tamu.edu

}

\end{document}